\documentclass{atmp}

\usepackage{url}
\usepackage[table]{xcolor}
\definecolor{lightgray}{gray}{0.9}
\usepackage{amssymb}

\newcommand{\be}{\begin{equation}} 
\newcommand{\ee}{\end{equation}}
\newcommand{\bea}{\begin{eqnarray}} 
\newcommand{\eea}{\end{eqnarray}}
\newcommand{\bc}{\begin{center}}
\newcommand{\ec}{\end{center}}

\newcommand{\fa}{\mathfrak{a}}
\newcommand{\fh}{\mathfrak{h}}
\newcommand{\g}{\mathfrak{g}}
\newcommand{\gl}{\mathfrak{gl}}
\newcommand{\C}{\mathbb{C}}

\begin{document}

\markboth{L.Fortunato and W.A. de Graaf}{E$_7 \subset $ Sp(56,${\mathbb R}$) irrep decomposition}

\title{E$_7 \subset $ Sp(56,${\mathbb R}$) irrep decompositions of interest for physical models}

\author{L.Fortunato\footnote{Dip. Fisica e Astronomia ``G.Galilei'', Univ. Padova and INFN-Sezione di Padova, via Marzolo 8, I-35131 Padova, Italy}, W.A. de Graaf\footnote{Dip. Matematica, Univ. Trento, via Sommarive 24, I-38123  Povo (Trento), Italy}}

\maketitle

\begin{abstract}
In this note we show how to obtain the projection matrix for the $E_7 \subset C_{28}$ chain and we tabulate some decompositions of the symplectic algebra $C_{28}$ representations into irreps of the $E_7$ subalgebra that are important for various physical models.
\end{abstract}
\maketitle

\section{Introduction and motivation}
The symplectic Lie groups and the corresponding algebras emerge in several physical models as symmetry group or dynamical symmetry groups. Among them, in particular, the group Sp(56,${\mathbb R}$), despite having an unusually  high rank, shows a distinctive trait, namely it contains $E_7$ as a subgroup. 
It appears, for example, in the study of symplectic extensions of 
$N=8$ gauged supergravity theories 
\cite{Dall}, 
that are extensions of the SO(8) gauged supergravity, as the group of transformations that act on the electric and magnetic gauge fields vectors.
These theories in particular need to look at the way in which certain representations branch into irreps of smaller algebras of distinguished importance and this fact motivated our effort to look for this particular branching that, due to its peculiarity, does not conform itself to known algorithms and requires a computational approach. 
  
High-rank symplectic algebras appear also prominently in various models of nuclear and atomic physics (cfr. for instance Ref. \cite{Geor}) and in many-body quantum theories \cite{Fort}. The Sp(2n,${\mathbb R}$) group is the dynamical group of the n-dimensional isotropic harmonic oscillator, that has a fundamental importance in all realms of quantum physics\cite{Wyb}.
While we are not interested here in the physical models based on these groups, nor what is their precise physical interpretation, they rely on the classification of the corresponding Lie subalgebras, a very interesting mathematical topic by itself (see, for example, Ref. \cite{DeGr}), and on the solution of the branching problem, that is very challenging and might become computationally demanding with increasing rank of involved algebras. 

We will investigate here the corresponding Lie algebras.
The exceptional Lie algebra $E_7$ is a special subalgebra of the $C_{28}$ symplectic algebra, i.e. $C_{28} \supset E_7$.
In order to study physical models based on this algebra chain, it is of fundamental importance to be able to determine 
the branching problem from the higher to the lower algebra. 

Notice that this particular branching rule does not appear in Ref. \cite{Laro}, where many other important branching rules are derived.
To the best of our knowledge the first treatment of the present subalgebra chain is given in the Ph.D. thesis of M. Lorente \cite{Lore}.

\section{Projection matrix}
The problem of finding the branching rules amounts to finding the projection matrix. 
In order to get this matrix we first describe how we realize the embedding $E_7\subset C_{28}$.  
We start with the Lie algebra of
type $E_7$, construct its lowest dimensional representation, and then build the Lie
algebra of type $C_{28}$ around that.

First we recall the following fact. Let $M_0$ be a $56\times 56$-matrix over $\C$,
such that $\det(M_0)\neq 0$ and $M_0^T = -M_0$. Set $\g = \{ a \in \gl(56,\C) \mid
a^TM_0 = -M_0 a \}$. Then $\g$ is a simple Lie algebra of type $C_{28}$. 

Let $\fa$ be the simple Lie algebra of type $E_7$ over $\C$. The computer algebra
system {\sf GAP}4 (\cite{gap}) contains a function for constructing the irreducible
representations of $\fa$. A description of the underlying algorithm can be found in Ref. 
\cite{wdg1}. We used it to compute the $56$-dimensional representation
$\rho : \fa \to \gl(56,\C)$. Next we considered the space $\mathcal{M}$ consisting of 
all $56\times 56$-matrices $M$ such that $M^T = -M$ and $\rho(x)^TM = -M\rho(x)$,
where $x$ runs through a basis of $\fa$. Note that this yields a set of linear
equations for the entries of $M\in \mathcal{M}$. From the solution of this system of equations,
 it turned out that $\dim \mathcal{M}=1$ and we let $M_0$ be a basis element of $\mathcal{M}$. 
It also turns out that $\det(M_0)\neq 0$.
So we constructed $\g$ as above, which is the simple Lie algebra of type $C_{28}$
containing $\rho(\fa)$. 

In order to describe the construction of the projection matrix we recall the following.
Let $\fh$ be a semisimple Lie algebra over $\C$ of rank $n$. Let $C$ be the Cartan
matrix of its root system. Then there are elements $x_i, y_i, h_i \in \fh$ such that
\begin{eqnarray*}
&&[h_i,h_j]=0\\ 
&&[x_i,y_j] = \delta_{i,j} h_i\\
&&[h_j,x_i] = C(i,j) x_i\\
&&[h_j,y_i] = -C(i,j)y_i.
\end{eqnarray*}
These elements generate $\fh$ and are said to form a {\em canonical generating set}
of $\fh$. 

We enumerate the nodes of the Dynkin diagrams of the Lie algebras $\fa$ and
$\g$ as follows:
 
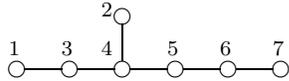
\begin{picture}(140,40)
  \put(23,0){\circle{6}}
\put(20,5){\footnotesize{1}}
  \put(43,0){\circle{6}}
\put(40,5){\footnotesize{3}}
  \put(63,0){\circle{6}}
\put(55,5){\footnotesize{4}}
  \put(83,0){\circle{6}}
\put(80,5){\footnotesize{5}}
  \put(63,20){\circle{6}}
\put(55,20){\footnotesize{2}}
  \put(26,0){\line(1,0){14}}
  \put(46,0){\line(1,0){14}}
  \put(66,0){\line(1,0){14}}
  \put(63,3){\line(0,1){14}}
\put(103,0){\circle{6}}
\put(100,5){\footnotesize{6}}
\put(86,0){\line(1,0){14}}
\put(123,0){\circle{6}}
\put(120,5){\footnotesize{7}}
\put(106,0){\line(1,0){14}}
\end{picture}

\begin{picture}(240,40)
  \put(23,0){\circle{6}}
\put(20,5){\footnotesize{1}}
  \put(43,0){\circle{6}}
\put(40,5){\footnotesize{2}}
  \put(63,0){\circle{6}}
\put(60,5){\footnotesize{3}}
  \put(26,0){\line(1,0){14}}
  \put(46,0){\line(1,0){14}}
  \put(66,0){\line(1,0){14}}
\put(90,0){$\ldots$}
\put(116,0){\line(1,0){14}}
\put(133,0){\circle{6}}
\put(130,5){\footnotesize{27}}
\put(158,0){\circle{6}}
\put(155,5){\footnotesize{28}}
\put(135,-2){\line(1,0){21}}
\put(135,2){\line(1,0){21}}
\put(141,-3){$<$}
\end{picture}
~\\

Correspondingly, we get canonical generating sets $x_i^\fa,y_i^\fa,h_i^\fa$,
$1\leq i\leq 7$ of $\fa$ and $x_j^\g,y_j^\g,h_j^\g$,
$1\leq i\leq 28$ of $\g$. Here the $h_i^\fa$, $h_j^\g$ are constructed so that 
each $h_i^\fa$ lies in the space spanned by the $h_j^\g$. In other words, there
are $a_{i,j}$ such that $h_i^\fa = \sum_{j=1}^{28} a_{i,j} h_j^\g$. Now the 
matrix $(a_{i,j})$ is called the projection matrix, and it is displayed below:

\be
\left(
\begin{array}{c}
 0 0 0 0 0 0 1 0 1 0 1 0 1 1 0 1 1 0 1 1 1 0 0 0 0 0 0 0  \\
 0 0 0 0 1 0 1 2 1 0 0 0 0 0 0 0 0 0 0 1 1 1 0 1 1 0 1 2  \\
 0 0 0 0 1 2 1 1 0 0 0 1 1 0 1 1 0 1 1 1 0 1 1 1 0 0 0 0  \\
 0 0 0 1 0 0 0 0 1 2 1 0 0 0 0 0 1 1 1 0 1 1 1 0 1 1 1 0  \\
 0 0 1 0 0 0 0 0 0 0 1 2 1 2 1 1 0 0 0 0 0 0 0 1 1 2 1 2  \\
 0 1 0 0 0 0 0 0 0 0 0 0 1 1 2 1 2 1 1 2 1 1 2 1 1 0 0 0  \\
 1 0 0 0 0 0 0 0 0 0 0 0 0 0 0 1 1 2 1 1 2 1 1 2 1 2 3 2 \\
 \end{array}
 \right)
\ee

\section{Branching rules}
This projection matrix is one of the main inputs to the algorithm for computing
branching rules, see \cite{navpat} and also \cite{wdg2}, Section 8.13. Given this
matrix the algorithm then uses the combinatorics of weights and root systems to 
obtain the decomposition of a given irreducible representation of $\g$,
when restricted to $\fa$. This algorithm has been implemented in the {\sf SLA}
package (\cite{sla}) inside the computer algebra system {\sf GAP}4 (\cite{gap}).
Using this we have investigated the branching problem for the smallest representations.

\begin{table}[t]
\caption{Branching $C_{28}\rightarrow E_7$ \label{tab:Irreps} for irrep of dimension smaller than one million.} 
\begin{center}
\rowcolors{1}{}{lightgray}
\begin{tabular}{ll|ll}
h.w. & Dim. & Branching & Dimensions \\ \hline
$[1,0_{27}]$ & {\bf 56} & $[0_6,1]$ & {\bf 56 }\\
$[0,1,0_{26}]$ &{\bf 1539} &$[0_5,1,0]$ & {\bf 1539} \\  
$[2,0_{27}]$ &{\bf 1596}    &   $[0_6,2] + [1,0_6]$&{\bf 1463}  +{\bf 133}\\
$[0_2,1,0_{25}]$ & {\bf 27664}  &$[0_4,1,0_2]$& {\bf 27664} \\
$[3,0_{27}]$ &{\bf 30856 }   &  $[0_6,3] + [1,0_5,1]+[0_6,1]$& {\bf 24320}  +{\bf 6480} +{\bf 56}\\
$[1,1,0_{26}]$ & {\bf 58464 }& $[ 0_5, 1, 1 ]+ [ 1, 0_5, 1 ]+[ 0,1, 0_5 ]$&{\bf 51072}+{\bf 6480}+{\bf 912}\\
$[0_3,1,0_{24}]$ &{\bf 365750} & $[ 0_3, 1, 0_3]$& {\bf 365750} \\
$[0,2,0_{26}]$ &{\bf 817740 } &  $[ 0_5, 2, 0 ]+ [ 1, 0_5, 2 ]+$ & {\bf 617253}+{\bf 150822} +\\
&& $[ 0, 1, 0_4, 1 ]+[ 2, 0_6 ]+[ 0_5, 1, 0 ]$ & {\bf 40755}+{\bf 7371}+{\bf 1539 }\\
\hline
\end{tabular}
\end{center}
\end{table}

We enumerate in table \ref{tab:Irreps} the branching of all the irreps of $C_{28}$ of dimension smaller than one million. We give the $C_{28}$ highest weight vectors in the first column, the dimension of the corresponding representation in the second column, the branching into $E_7$ representations in the third column and the corresponding dimensions in the fourth column. All these representations appear with multiplicity one. The notation $0_k$ in the highest weight vectors simply stands for $\cdots,0,\cdots$ repeated $k$ times.
We find, for example, that the 1596-dimensional adjoint representation of $C_{28}$ is decomposed into the sum of the 1463-dimensional representation and 133-dimensional adjoint representation of $E_7$, a most natural finding that excludes, however, the possibility to have other partitions of 1596 into the set $\{1,56,133,912,1463,1539\}$. In general, one might get an idea of the combinatorial complexity of the problem from the solution of the Frobenius Diophantine equation: excluding singlets, one gets three such possibilities (namely {\bf 1596}=12$\times${\bf 133}, 19$\times${\bf 56}+4$\times${\bf 133} and {\bf 1463}+{\bf 133}) and the number of solutions raises to 240 possible partitions when including the one-dimensional representation (i.e. $[0_{28}]$). Of course not all of these partitions are admissible from a purely algebraic perspective and the branching rules indicate how to accomplish the task with the correct algorithm that embeds representations of the smaller algebra into representations of the larger one. Due to the dimensions of the algebras involved, these numbers might grow very rapidly for higher representations.
Notice that, among the 27 possible representations of $E_7$ that have smaller
dimension than 1 million, only 14 appear in the decomposition.
(A more complete list of these representations can be found in the tables of Ref. \cite{Feg}).

It is also interesting to know how tensor products of the $C_{28}$ algebra representations decompose into sum of $C_{28}$ irreps: for example the {\bf 56}$\times${\bf 56} (3136-dimensional) representation decomposes into the sum of {\bf 1596}+{\bf 1539}+{\bf 1}. This kind of decompositions are easily accomplished with {\sf GAP} when the projection matrix is known. 
We enumerate the smallest in Table \ref{tab:Tens} because they might be useful to the investigation of physical models based on the $E_7 \subset C_{28}$ subalgebra chain.

\begin{table}[t]
\caption{\label{tab:Tens} 
Decomposition of tensor products of $C_{28}$ irreps of dimension smaller than one million. We have used here the dimensional notation, coefficients in the last two lines stand for multiplicities.} 
\begin{center}
\rowcolors{1}{}{lightgray}
\begin{tabular}{l|c|l}
Tensor Product & Dimension & Decomposition \\ \hline
{\bf 56}$\times${\bf 56}  & 3136&{\bf 1596 }+ {\bf 1539}+ {\bf 1}\\
{\bf 56}$\times${\bf 1539}& 86184& {\bf 58464}+{\bf 27664 }+  {\bf 56}\\
{\bf 56}$\times${\bf 1596}& 89376& {\bf 58464}+{\bf 30856 }+  {\bf 56}\\
{\bf 56}$\times${\bf 56}$\times${\bf 56}  & 175616 &  2({\bf 58464})+ {\bf 30856}+ \\
&&+{\bf 27664}+3({\bf 56})\\
\hline
\end{tabular}
\end{center}
\end{table}

\section*{Acknowledgments}
We thank G.Dall'Agata (Padova) for having called our attention to this problem and for useful discussions.


\begin{thebibliography}{zzz}
\bibitem{Dall} G. Dall'Agata, G. Inverso, M. Trigiante, Phys. Rev. Lett. {\bf 109}, 201301 (2012)
\bibitem{Geor} A. I. Georgieva, et al. J. Phys. Conf. Ser. 284, 012028 (2011)
\bibitem{Fort} L. Fortunato, AIP Conf. Proc. 1488, 366 (2012)
\bibitem{Wyb} B. G. Wybourne, {\it Classical Groups for Physicists}, John Wiley \& sons Inc., New York, (1974)
\bibitem{Laro} M. Larouche, J. Patera, J.Phys. A: Math.Theor., 44, 115203 (2011) 
\bibitem{Lore} M. Lorente-Paramo, Ph.D. thesis, Un. Compl. Madrid (1972), in spanish.
\bibitem{DeGr} W. A. de Graaf, J. Algebra, {\bf 325} 416-430 (2011)
\bibitem{wdg1} W. A. de Graaf, J. Pure And Applied Algebra {\bf 164}, 87-107 (2001)
\bibitem{wdg2} W. A. de Graaf, {\it Lie Algebras: Theory and Algorithms}, Elsevier, Amsterdam (2000)
\bibitem{navpat} A. Navon, J. Patera, J. Math. Phys. {\bf 8}, 489-493 (1967)
\bibitem{gap}  The GAP Group, GAP -- Groups, Algorithms, and Programming, Version 4.5; \url{http://www.gap-system.org}
\bibitem{sla} W. A. de Graaf, {\it SLA - computing with Simple Lie Algebras. A GAP package}, \url{http://science.unitn.it/~degraaf/sla.html} 
\bibitem{Feg} R. Feger, T.W. Kephart, arXiv:1206.6379 (2012)
\end{thebibliography}
\end{document}